\documentclass[a4paper,10pt]{article}

\usepackage[latin1]{inputenc}		
\usepackage{amsfonts}			
\usepackage{amsmath}			
\usepackage{amsthm}			
\usepackage{amssymb}			
\usepackage{graphicx}			
\usepackage{enumerate}			
\usepackage{mathrsfs}			

\newcommand{\nN}{\mathbb{N}}
\newcommand{\nZ}{\mathbb{Z}}
\newcommand{\nQ}{\mathbb{Q}}
\newcommand{\nR}{\mathbb{R}}
\newcommand{\nC}{\mathbb{C}}
\newcommand{\nD}{\mathbb{D}}
\newcommand{\nP}{\mathbb{P}}

\newcommand{\lhs}{\mathbb{I}}
\newcommand{\rhs}{\mathbb{II}}

\theoremstyle{plain}
\newtheorem{mythm}{Theorem}[section]
\newtheorem{mylem}[mythm]{Lemma}
\newtheorem{myprop}[mythm]{Proposition}
\newtheorem{mycor}[mythm]{Corollary}

\theoremstyle{definition}
\newtheorem{mydef}[mythm]{Definition}
\newtheorem{myoss}[mythm]{Remark}
\newtheorem{myes}[mythm]{Example}

\newtheorem{myces}[mythm]{Counterexample}

\theoremstyle{remark}

\newenvironment{mydim}[0]		
{\begin{proof}}
{\end{proof}}

\newenvironment{mytext}[0]{}{}		

\newcommand{\subrelcpt}{{\subset\subset}}

\newcommand{\mc}[1]{\mathcal{#1}}

\newcommand{\mg}[1]{\mathfrak{#1}}

\newcommand{\on}[1]{\operatorname{#1}}

\newcommand{\sumaa}[2]{\hspace{-#2}\sum_{#1}\hspace{-#2}}

\newcommand{\abs}[1]{\left|#1\right|}

\newcommand{\norm}[1]{\left\|#1\right\|}

\newcommand{\vect}[1]{\overrightarrow{#1}}

\newenvironment{myabstract}[1]
{\begin{center}\small{\textbf{#1}}

\vspace{0.25cm}
\begin{minipage}[l]{10.4cm}\begin{small}}
{\end{small}\end{minipage}\end{center}}

\newcommand{\mysubsubsect}[1]
{
\subsubsection*{#1}
}

\newcommand{\mysubsect}[1]
{
\subsection*{#1}
\addcontentsline{toc}{subsection}{#1} 
}

\newcommand{\mysect}[1]
{
\section{#1}
\markright{\hfill \thesection\ #1 \hfill}
}

\newcommand{\mysecta}[1]
{
\section*{#1}
\markright{\hfill \thesection\ #1 \hfill}
}

\newcommand{\mydefal}[2]{\textbf{#1}} 
\newcommand{\mydefin}[1]{\textbf{#1}}

\title{Rigidification of holomorphic germs with non-invertible differential}
\author{Matteo Ruggiero\textsuperscript{a,}\footnote{Corresponding author. Tel: (+39) 050 509092. Fax: (+39) 050 563513.} \textsuperscript{,}\footnote{E-mail: m.ruggiero@sns.it}
\\ \footnotesize{\textsuperscript{a} Scuola Normale Superiore, Piazza dei Cavalieri 7, 56126 Pisa, Italy.}}
\date{}

\makeindex

\begin{document}

\maketitle

\begin{myabstract}{Abstract}
We study holomorphic germs $f:(\nC^2, 0) \rightarrow (\nC^2,0)$ with non-invertible differential $df_0$.
In order to do this, we search for a modification $\pi:X \rightarrow (\nC^2,0)$ (i.e., a composition of point blow-ups over the origin), and an infinitely near point $p \in \pi^{-1}(0)$, such that the germ $f$ lifts to a holomorphic germ $\hat{f}:(X,p) \rightarrow (X,p)$ which is rigid (i.e., the generalized critical set of $\hat{f}$ is totally invariant and has normal crossings at $p$).
We extend a previous result (see \cite{favre-jonsson:eigenval}) for superattracting germs to the general case, and deal with the uniqueness of this process in the semi-superattracting case ($\on{Spec}(df_0)=\{0, \lambda\}$ with $\lambda \neq 0$).
We specify holomorphic normal forms for the nilpotent case and for the type $(0,\nD)$, that is $\on{Spec}(df_0)=\{0, \lambda\}$ with $\lambda$ in the unitary disk $\nD \subset \nC$, and formal normal forms for the type $(0, \nC \setminus \nD)$.
\end{myabstract}

\mysecta{Introduction}

Our aim in this paper is to study the structure of non-invertible holomorphic germs $f:(\nC^2, 0) \rightarrow (\nC^2,0)$.

We shall consider only dominant holomorphic germs:

\begin{mydef}
Let $f:(\nC^2, 0) \rightarrow (\nC^2,0)$ be a holomorphic germ.
Then $f$ is \mydefin{dominant} if $\on{det}(df_p)$ is not identically zero. 
\end{mydef}

In this paper we are particularly interested in the following classes of holomorphic germs:
\begin{mydef}
Let $f:(\nC^2, 0) \rightarrow (\nC^2,0)$ be a holomorphic germ, and let us denote by $\on{Spec}(df_0)=\{\lambda_1, \lambda_2\}$ the set of eigenvalues of $df_0$.
Then $f$ is said:
\begin{itemize}
\item \mydefin{attracting} if $\abs{\lambda_i} < 1$ for $i=1, 2$;
\item \mydefin{superattracting} if $df_0=0$;
\item \mydefin{nilpotent} if $df_0$ is nilpotent (i.e., $df_0^2=0$; in particular, superattracting germs are nilpotent germs);
\item \mydefin{semi-superattracting} if $\on{Spec}(df_0)=\{0, \lambda\}$, with $\lambda \neq 0$;
\item \mydefin{of type $(0,D)$} if $\on{Spec}(df_0)=\{0, \lambda\}$ and $\lambda \in D$, with $D \subset \nC$ a subset of the complex plane.
\end{itemize}
In particular the semi-superattracting germs are the ones of type $(0,\nC^*)$.
\end{mydef}
We shall denote by $\nD$ the open disk of radius $1$ centered at $0$.

A typical problem one would like to solve is to find a classification up to local (holomorphic, formal or topological) conjugacy; while this problem is mostly solved in dimension $1$, there are classifications of germs in dimension $2$ only in a few cases. One of these cases is the formal and holomorphic classification of attracting rigid germs, proved by C. Favre (see \cite{favre:rigidgerms}).

\begin{mydef}
Let $f:(\nC^n, 0) \rightarrow (\nC^n,0)$ be a (dominant) holomorphic germ.
We denote by $\mc{C}(f)=\{z\ |\ \det (df_z)=0\}$ the \mydefal{critical set}{critical set} of $f$, and by $\mc{C}^\infty(f)=\bigcup_{n\in \nN} f^{-n}\mc{C}(f)$ the \mydefal{generalized critical set}{generalized critical set} of $f$.
Then a (dominant) holomorphic germ $f$ is \mydefal{rigid}{rigid germ} if:
\begin{enumerate}[(i)]
\item $\mc{C}^\infty(f)$ (is empty or) has normal crossings at the origin; and
\item $\mc{C}^\infty(f)$ is forward $f$-invariant.
\end{enumerate}
\end{mydef}

\begin{myoss}
In \cite{favre:rigidgerms}, the condition (ii) is not explicitly stated in the definition of a rigid germ, but it is implicitly used.
The second property does not follow from the first one: if for example we consider the map $f(z,w)=(\lambda z^p, z(1+w^2))$, with $p \geq 1$ and $\lambda \in \nC^*$, the generalized critical set is $\{zw=0\}$, but $f(z,0)=(\lambda z^p, z)$ and hence $\mc{C}^\infty(f)$ is not forward $f$-invariant.
\end{myoss}

A possible way for studying the local dynamics of a generic holomorphic germ in $(\nC^2,0)$, and for finding some invariants up to conjugacy, is suggested by the continuous local dynamics (see \cite[Chapters 1 and 2]{ilyashenko-yakovenko:analDE} for main techniques in continuous local dynamics, and \cite{seidenberg:redsing} for Seidenberg's Theorem): we can blow-up the fixed point (the origin), replacing the ambient space by a more complicated space but simplifying the map, and study the lift $\hat{f}$ of $f$. But a single blow-up is often not enough, and one is led to consider a composition of point blow-ups $\pi:X \rightarrow (\nC^2,0)$ over the origin (called \mydefin{modification}).

A clever way to study all modifications at the same time was introduced by Favre and Jonsson (see \cite{favre-jonsson:valtree}). Let us take the set of all modifications $\mg{B}$; for every $\pi \in \mg{B}$ we can consider a simplicial graph $\Gamma_\pi^*$, whose vertices are the irreducible components of the exceptional divisor of $\pi$ (we shall call these vertices \mydefin{exceptional components}).
Taking the direct limit of these simplicial graphs, we obtain the ($\nQ$-)\mydefin{universal dual graph} $\Gamma^*$, which has a natural $\nQ$-tree structure.
Since it is easier to work with $\nR$-trees, we can take the completion $\Gamma$ of $\Gamma^*$, called the \mydefin{universal dual graph}.

Favre and Jonsson also showed that the universal dual graph is strictly related to the set $\mc{V}$ of all centered and normalized valuations on the ring of formal power series in $2$ coordinates: $\mc{V}$ admits an $\nR$-tree structure, and is isomorphic to $\Gamma$ (in a strong sense).

Basically thanks to this isomorphism, that relates the geometry of exceptional components to the algebra of valuations, it is possible to define the action $f_\bullet: \mc{V} \rightarrow \mc{V}$ on the valuative tree $\mc{V}$ induced by a holomorphic germ $f:(\nC^2, 0)\rightarrow (\nC^2,0)$ (see \cite{favre-jonsson:eigenval}).

Favre and Jonsson studied in \cite{favre-jonsson:eigenval} the dynamical behavior of $f_\bullet$ when $f:(\nC^2, 0) \rightarrow (\nC^2,0)$ is superattracting; in particular, they proved that one can find a modification $\pi:X \rightarrow (\nC^2, 0)$ and a point $p \in \pi^{-1}(0)$ such that the lift $\hat{f}:(X,p) \rightarrow (X,p)$ defined as a birational map by $\hat{f}=\pi^{-1}\circ f \circ \pi$ is actually holomorphic in $p$, and rigid.
This can be done by finding a fixed point $\nu_\star$ for $f_\bullet$, called eigenvaluation, and studying the basin of attraction around this eigenvaluation.
We shall call this process \mydefin{rigidification}.

\begin{mydef}
Let $f:(\nC^2,0) \rightarrow (\nC^2,0)$ be a (dominant) holomorphic germ. Let $\pi:X \rightarrow (\nC^2,0)$ be a modification and $p \in \pi^{-1}(0)$ a point in the exceptional divisor of $\pi$.
Then we shall call the triple $(\pi,p,\hat{f})$ a \mydefin{rigidification} of $f$ if the lift $\hat{f}=\pi^{-1} \circ f \circ \pi$ is a holomorphic rigid germ in $p$.
\end{mydef}

We shall follow their strategy for finding eigenvaluations and rigidifications, extending the result to all (dominant) holomorphic germs.
We notice that the rigidification process is trivial if $df_0$ is invertible (because the map $f$ is itself rigid).
Our main result can be stated as:

\begin{mythm}\label{thm:rigidification}
Every (dominant) holomorphic germ $f:(\nC^2,0)\rightarrow (\nC^2,0)$ admits a rigidification.
\end{mythm}

The nilpotent case is quite the same as the superattracting case dealt in \cite{favre-jonsson:valtree} (see Remark \ref{oss:superattrvsnilpot}). Hence we shall focus on the semi-superattracting case, proving a sort of uniqueness of the rigidification process in this case, that can be stated as follows.

\begin{mythm}\label{thm:fixedvaluniqueness}
Let $f$ be a (dominant) semi-superattracting holomorphic germ. Then $f$ admits a unique eigenvaluation $\nu_\star$, that has to be a (possibly formal) curve valuation with multiplicity $m(\nu_\star)=1$. Let us denote $\nu_\star=\nu_C$, with $m(C)=1$. Then, (only) one of the following holds:
\begin{enumerate}[(i)]
\item the set of valuations fixed by $f_\bullet$ consists only of the eigenvaluation $\nu_\star$, there exists (only) one contracted critical curve valuation $\nu_D$, and in this case, it has to be $m(D)=1$;
\item the set of valuations fixed by $f_\bullet$ consists of two valuations, the eigenvaluation $\nu_\star$, and a curve valuation $\nu_D$, where $D$ is a (possibly formal) curve with $m(D)=1$.
\end{enumerate}
In both cases, $C$ and $D$ have transverse intersection, i.e., their intersection number is $C \cdot D=1$.
\end{mythm}

We shall then prove the formal classification of semi-superattracting rigid germs (the first case actually follows from the holomoprhic classification of such germs given in \cite{favre:rigidgerms}).

\begin{mythm}\label{thm:rigidsemiattr}
Let $f:(\nC^2, 0) \rightarrow (\nC^2,0)$ be a (holomorphic) semi-superattracting rigid germ. Let $\lambda \in \nC^*$ be the non-zero eigenvalue of $df_0$.
\begin{enumerate}[(i)]
\item If $\abs{\lambda} < 1$ or $\lambda= e^{2 \pi i \theta}$ with $\theta \in \nR \setminus \nQ$, then $f$ is formally conjugated to the map
$$
(z,w) \mapsto (\lambda z, z^c w^d)\mbox.
$$
\item If $\abs{\lambda} > 1$, then $f$ is formally conjugated to the map
$$
(z,w) \mapsto \big(\lambda z, z^c w^d(1+\varepsilon z^l)\big)\mbox,
$$ 
where $\varepsilon \in \{0,1\}$ if $\lambda^l=d$ (\mydefin{resonant case}), and $\varepsilon=0$ otherwise.
\item If there exists $r \in \nN^*$ such that $\lambda^r=1$, then $f$ is formally conjugated to the map
$$
(z,w) \mapsto \big(\lambda z (1 + z^s + \beta z^{2s}), z^c w^d(1+\varepsilon(z^r))\big)\mbox,
$$
where $r|s$ and $\beta \in \nC$, while $\varepsilon$ is a formal power series in $z^r$, and $\varepsilon \equiv 0$ if $d \geq 2$.
\end{enumerate}
In all cases, $c \geq 0$, $d \geq 1$ and $c+d \geq 2$.
\end{mythm}

\begin{myoss}
In Theorem \ref{thm:rigidsemiattr}.(ii), understanding in the resonant case which one of the two possible normal forms (with $\varepsilon = 0$ or $1$) is the normal form of a given germ $f$ seems quite difficult, either by the dynamics of $f$, or by the action of $f_\bullet$ (see also Remark \ref{oss:resonancebulletaction}).
\end{myoss}

We shall also present two counterexamples (see Counterexamples \ref{ces:nonhol1} and \ref{ces:nonhol2}) that show how the holomorphic classification of rigid germs of type $(0, \nC \setminus \nD)$ is not trivial, meaning that it does not coincide with the formal classification.

Finally we shall use the holomorphic classification of attracting rigid germs given in \cite{favre:rigidgerms} to give holomorphic normal forms for a rigidification for type $(0,\nD^*)$ (see Proposition \ref{prop:rigidattr}), and Theorem \ref{thm:rigidsemiattr} to give formal normal forms for a rigidification for type $(0,\nC \setminus \nD)$ (see Proposition \ref{prop:rigidrepel} and Proposition \ref{prop:rigidindif}).

\begin{mytext}
Using a different language, Theorem \ref{thm:rigidification} says that one can suppose a germ to be rigid, up to birational conjugacy. Then the normal forms of a rigidification give us normal forms for the birational classification of these germs.
In the semi-superattracting case, we prove Theorem \ref{thm:fixedvaluniqueness}, a sort of uniqueness of this process, that leads to a sort (see Example \ref{es:normalformchanges}) of uniqueness for these normal forms.
The dynamics of these rigidifications $\hat{f}$, easier to study than the initial germ $f$ itself, give us informations on the dynamics of $f$ (by projection), while the birational classification gives us informations on the holomorphic classification, in a very consistent way.
In fact (see Remark \ref{oss:uniqueeigenval}) the action of $\hat{f}_\bullet$ is related to the action of $f_\bullet$ in a suitable basin of attraction in the valuative tree.
\end{mytext} 

\begin{mytext}
This paper is divided into four sections. In the first section we recall the construction of the valuative tree $\mc{V}$ and its isomorphic equivalent, the universal dual graph $\Gamma$, as in \cite{favre-jonsson:valtree}; the action $f_\bullet$ induced by a (dominant) holomorphic germ $f$; the existence of an eigenvaluation and of a basin of attraction, as in \cite{favre-jonsson:eigenval}, adapted to deal with the general case.
In the second section we prove Theorem \ref{thm:rigidification} and Theorem \ref{thm:fixedvaluniqueness}.
In the third section we deal with the classification of semi-superattracting rigid germs.
In the last section we compute normal forms for a rigidification in every case, and we end with some remarks on the rigidification process.
\end{mytext}

\mysubsect{Acknowledgements}
The author would like to thank Charles Favre and Mattias Jonsson for several useful remarks and clarifications on their works, and Marco Abate for many suggestions.
 
\mysect{The Valuative Tree}

\mysubsect{Modifications}

\begin{mytext}
The main objects we want to study are modifications, i.e., compositions of point blow-ups, and the lifts of maps over the exceptional divisor of a modification.
Hence let us start fixing notations.
\end{mytext}

\begin{mydef}
Let $X$ be a complex $2$-manifold, and $p \in X$ a point. We call a holomorphic map $\pi:Y \rightarrow (X,p)$ a \mydefal{modification}{modification} over $p$ if $\pi$ is a composition of point blow-ups, with the first one being over $p$, and such that $\pi$ is a biholomorphism outside $\pi^{-1}(p)$.
We call $\pi^{-1}(p)$ the \mydefal{exceptional divisor}{exceptional divisor} of $\pi$, and we call every irreducible component of the exceptional divisor an \mydefal{exceptional component}{exceptional component}.
We will denote by $\mg{B}$ the set of all modifications over $0 \in \nC^2$, and by $\Gamma_\pi^*$ the set of all exceptional components of a modification $\pi$.
We will call a point $p \in \pi^{-1}(0)$ on the exceptional divisor of a modification $\pi \in \mg{B}$ an \mydefin{infinitely near point} (we consider $0 \in \nC^2$ as an infinitely near point too).
\end{mydef}

\mysubsect{Tree Structure}

\begin{mytext}
We shall here fix notations for $\nR$-trees; see \cite[Chapter 3]{favre-jonsson:valtree} for definitions and proofs, and \cite[Section 4]{favre-jonsson:eigenval} for properties of tree maps.
\end{mytext}

\begin{mydef}
Let $(\mc{T}, \leq)$ be an $\nR$-tree. 
Maximal elements of $\mc{T}$ shall be called \mydefal{ends}{tree!ends}.

Let $\tau_1, \tau_2 \in \mc{T}$ be two points. We shall denote by $[\tau_1,\tau_2]$ (resp., $[\tau_1, \tau_2)$ and $(\tau_1, \tau_2)$) the \mydefin{closed} (resp., \mydefin{semiopen} and \mydefin{open}) \mydefal{segment}{tree!segment} between $\tau_1$ and $\tau_2$.

We shall denote by $T_\tau \mc{T}$ the \mydefal{tangent space}{tree!tangent space} of $\mc{T}$ over a point $\tau$, and we denote by $\vect{v}=[\sigma] \in T_\tau \mc{T}$ a tangent vector over $\tau$ (represented by $\sigma$).
Then the point $\tau$ is a \mydefal{terminal point}{tree!terminal point}, a \mydefal{regular point}{tree!regular point} or a \mydefal{branch point}{tree!branch point} if $T_\tau \mc{T}$ has $1$, $2$ or more than $2$ tangent vectors respectively.

Finally, let $\tau \in \mc{T}$ be a point in the tree, and $\vect{v}=T_\tau \mc{T}$ a tangent vector over it; we shall denote by
$$
U_\tau(\vect{v}):=\{\sigma \in \mc{T}\ |\ \vect{v}=[\sigma]\}\mbox{.}
$$
the \mydefin{(weakly) open set} associated to $\vect{v}$ in $\tau$.
\end{mydef}

\mysubsect{Universal Dual Graph}

\mysubsubsect{Dual Graph of a Modification}

\begin{mytext}
Given a modification $\pi \in \mg{B}$, we can equip the set $\Gamma_\pi^*$ of all exceptional components of $\pi$ with a simplicial tree structure (i.e., an $\nN$-tree structure, see \cite[pages 51, 52]{favre-jonsson:valtree}).
\end{mytext}

\begin{mydef}
We fix the set of vertices ($\Gamma_\pi^*$), and we say that two exceptional components are joined by an edge if and only if their intersection is non-empty.
We will denote by $\leq_\pi$ the induced partial ordering (given by the corrispondence between simplicial trees and $\nN$ trees).
Then $(\Gamma_\pi^*, \leq_\pi)$ will be called the \mydefin{dual graph} of $\pi$.
\end{mydef}

\begin{mydef}
Let $\pi \in \mg{B}$ be a modification. A point $p \in \pi^{-1}(0)$ in the exceptional divisor of $\pi$ is a \mydefin{free point} (respectively a \mydefin{satellite point}) if $\pi$ is a regular point (respectively a singular point) of $\pi^{-1}(0)$.
\end{mydef}

\begin{mytext}
We notice that satellite points are also known as \mydefin{corners} in literature.
An equivalent definition will be that $p$ is a free point if it belongs to only one exceptional component, while it is a satellite point if it belongs to exactly two exceptional components (that will have only one intersection point, with transverse intersection).
\end{mytext}

\mysubsubsect{Universal Dual Graph}

\begin{mydef}
We will call \mydefin{universal dual graph} the direct limit of dual graphs along all modifications in $\mg{B}$:
$$
(\Gamma^*, \leq):= \lim_{\stackrel{\longrightarrow}{\pi \in \mg{B}}} (\Gamma_\pi^*, \leq_\pi) \mbox{.}
$$
\end{mydef}

\begin{mytext}
The universal dual graph is a way to see all exceptional components of all the possible modifications, all together at the same time.
The next result follows from this construction.
\end{mytext}

\begin{myprop}[{\cite[Proposition 6.2, Proposition 6.3]{favre-jonsson:valtree}}]\label{prop:undualgraph}
The universal dual graph $\Gamma^*$ is a $\nQ$-tree, rooted at $E_0$ the exceptional component that arises from the single blow-up of the origin $0 \in \nC^2$. Moreover all points are branch points for $\Gamma^*$.
If we have an exceptional component $E \in \Gamma^*$, then $p \mapsto \vect{v_p}=[E_p]$, where $[E_p] \in T_E \Gamma^*$ is the tangent vector represented by the exceptional component that arises from the blow-up of $p$, gives a bijection from $E$ to $T_E \Gamma^*$.
\end{myprop}

\begin{mytext}
One can complete $\Gamma^*$ to a complete $\nR$-tree $\Gamma$, that will also be called the (complete) \mydefin{universal dual graph}.
\end{mytext}

\begin{mytext}
The (complete) universal dual graph is a very powerful tool, also thanks to all the structure that arises from the completeness of $\nR$. But we do not know how a holomorphic germ $f$ acts on the universal dual graph. The answer to this question can be given thanks to the algebraic equivalent to the universal dual graph, the \mydefin{valuative tree}.  
\end{mytext}

\mysubsect{Valuations}

\begin{mytext}
We shall denote by $R=\nC[[x,y]]$ the ring or formal power series in $2$ coordinates, and by $K=\nC((x,y))$ the quotient field of $R$, that is the field of Laurent series in $2$ coordinates.
Then $R$ is an UFD local ring, with maximal ideal $\mg{m}=\langle x,y \rangle$.
Favre and Jonsson considered a slightly different concept of valuation, that takes values in $[0, +\infty]$, while the classical Krull valuations take values into a (totally ordered) abelian group.
Moreover they focus their attention on \mydefin{centered} valuations, i.e., valuations $\nu:R \rightarrow [0, +\infty]$ that take strictly positive values on $\mg{m}$.
\end{mytext}

\begin{mytext}
The set of all (centered) valuations can be endowed by a partial order.
\end{mytext}

\begin{mydef}
Let $\nu_1, \nu_2$ be two centered valuations; then $\nu_1 \leq \nu_2$ if and only if $\nu_1(\phi) \leq \nu_2(\phi)$ for every $\phi \in R$.
\end{mydef}

\begin{mytext}
The set of all (normalized) centered valuations with this partial order admits an $\nR$-tree structure: we shall call $\mc{V}$ the \mydefal{valuative tree}{valuative tree}.
\end{mytext}

\begin{mytext}
Valuations are naturally embedded into Krull valuations, but the converse is not true (see the exceptional curve valuations, \cite[page 18]{favre-jonsson:valtree} for details).
\end{mytext}

\begin{mytext}
The next theorem is a classic result of algebraic geometry: see \cite[Part VI, Chapter 5]{zariski-samuel:commalg2}, or \cite{hartshorne:alggeo} for a modern exposition.
\end{mytext}
\begin{mythm}[{\cite[Theorem 4.7]{hartshorne:alggeo}}]
Let $\nu$ be a Krull valuation on $K=\nC((x,y)) \supset \nC(x,y)$, $R_\nu$ the associated valuation ring and $\pi:X \rightarrow (\nC^2,0)$ a modification.
Then there exists a unique irreducible submanifold $V$ of $X$ such that $R_\nu$ dominates $\mc{O}_{X,V}$ the ring of regular functions in $V$. Moreover if $\nu$ is centered, then $V$ is a point or an exceptional component in $\pi^{-1}(0)$.

This $V$ is called the \mydefin{center} of $\nu$ in $X$.
\end{mythm}

\begin{mytext}
The center of a valuation is the main concept allowing to pass from valuations to exceptional components, and gives an isomorphism  between the valuative tree and the universal dual graph (see \cite[Theorem 6.22]{favre-jonsson:valtree}).
Moreover, the center of a valuation gives special open sets in the valuative tree (see \cite[Corollary 6.34]{favre-jonsson:valtree} for some properties).
\end{mytext}

\begin{mydef}
Let $p \in \pi^{-1}(0)$ be an infinitely near point of a modification $\pi:X \rightarrow (\nC^2, 0)$; we shall denote by $U(p) \subseteq \mc{V}$ the (weakly open) set of all valuations whose center in $X$ is $p$.
\end{mydef}

\begin{mytext}
For proofs and further details on valuations, see \cite[Part VI]{zariski-samuel:commalg2}.
\end{mytext}

\mysubsect{Classification of Valuations}

\begin{mytext}
We shall describe the classification of valuations and their role in the valuative tree (see \cite[Chapter 1]{favre-jonsson:valtree}). \end{mytext}

\mysubsubsect{Divisorial Valuations}

\begin{mytext}
Divisorial valuations are associated to an exceptional component $E$ of a modification $\pi$; in particular $\nu_E$ is defined by
$$
\nu_E(\phi):= (1/b_E) \on{div}_E(\pi^* \phi)\mbox{,}
$$
where $\on{div}_E$ is the vanishing order along $E$, $\pi^*\phi=\phi \circ \pi$ and $1/b_E$ is necessary to have a normalized valuation ($b_E \in \nN^*$ is known as the \mydefin{generic multiplicity} of $\nu_E$, see \cite[page 64]{favre-jonsson:valtree}, or the second Faray weight of $E$, see \cite[page 122]{favre-jonsson:valtree}).
The set of all divisorial valuations is often denoted by $\mc{V}_{div}$.
The divisorial valuations are the branch points of the valuative tree, and in particular we have $T_{\nu_E} \mc{V}\cong E$.

The most important example is the \mydefal{multiplicity valuation}{multiplicity valuation}, defined by
$$
\nu_\mg{m}(\phi):= m(\phi) = \max\{n\ |\ \phi \in \mg{m}^n\}\mbox{;} 
$$
it is associated to a single blow-up over the origin, and plays the role of the root of $\mc{V}$.
We will write $\nu_\mg{m}$ if we want to consider the multiplicity as a valuation (or better, as a point on the valuative tree), and $m$ if we want to consider only the multiplicity of an element of $R=\nC[[x,y]]$.
\end{mytext}

\mysubsubsect{Irrational Valuations}

\begin{mytext}
Irrational valuations are the regular points of the valuative tree.
Divisorial and irrational valuations are called \mydefal{quasimonomial valuations}{quasimonomial valuations}, their set will be denoted by $\mc{V}_{qm}$. 
For a geometric interpretation of quasimonomial valuations, see \cite[pages 16, 17]{favre-jonsson:valtree}.

Important examples of quasimonomial valuations are \mydefin{monomial valuations}.
Fix local coordinates $(x,y)$; then the monomial valuation of weights $(s,t)$ is defined by
$$
\nu_{s,t}\left(\sum_{i,j}a_{i,j}x^i y^j\right)= \min\{si+tj\ |\ a_{i,j} \neq 0\}\mbox{.}
$$
\end{mytext}

\mysubsubsect{Curve Valuations}

\begin{mytext}
Curve valuations are ends of the valuative tree.
They are associated to a (formal) irreducible curve (germ) $C=\{\psi=0\}$; in particular $\nu_C$ is defined by
$$
\nu_C(\phi):= \frac{C \cdot \{\phi=0\}}{m(C)}\mbox{,}
$$
where with $C \cdot D$ we denote the standard intersection multiplicity between the curves $C$ and $D$, and $m(C)=m(\psi)$ is the multiplicity of $C$ (in $0$).
We will often use the notation $\nu_\psi$ instead of $\nu_C$.

Analytic and non-analytic curve valuations have the same algebraic behavior, but they will play a different role as eigenvaluations, as we shall see in the proof of Theorem \ref{thm:rigidification}.
\end{mytext}

\mysubsubsect{Infinitely Singular Valuations}

\begin{mytext}
Infinitely singular valuations are the ones with $\on{rk}\nu = \on{ratrk}\nu = 1$ and $\on{trdeg}\nu = 0$, and share with curve valuations the role of ends of the valuative tree. 

It is not so simple to give a geometric interpretation of infinitely singular valuations, but we can think them as curve valuations associated to ``curves'' of infinite multiplicity.
They can be recognized also as valuations with infinitely generated value groups.
\end{mytext}

\mysubsect{Parametrizations}

\begin{mytext}
The valuative tree admits (at least) two natural parametrizations (skewness and thinness) and a concept of multiplicity, very useful for example to distinguish the type of valuations. For definitions and properties we refer to \cite[Chapter 3]{favre-jonsson:valtree}; all we need for the paper is the following.
\end{mytext}

\begin{myprop}[{\cite[Theorem 3.46]{favre-jonsson:valtree}}]
The thinness $A: \mc{V} \rightarrow [2, \infty]$ is a parametrization for the valuative tree. Moreover:
\begin{enumerate}[(i)]
\item the multiplicity valuation is the only one with $A(\nu_\mg{m})=2$;
\item for divisorial valuations we have $A(\nu_E) \in \nQ$;
\item for irrational valuations we have $A(\nu) \in \nR \setminus \nQ$;
\item for curve valuations we have $A(\nu_C)=\infty$;
\item for infinitely singular valuations we have $A(\nu) \in (2, \infty]$.
\end{enumerate}
\end{myprop}

\mysubsect{Dynamics on the Valuative Tree}

\mysubsubsect{Definition}

\begin{mytext}
In this section we will define the action $f_\bullet:\mc{V} \rightarrow \mc{V}$ induced by a holomorphic germ $f:(X,p) \rightarrow (Y,q)$ (where $X$ and $Y$ are two complex $2$-manifolds).
We shall also assume that $f$ is \mydefin{dominant}, i.e., $\on{rk} df$ is not identically $\leq 1$ near $p$.
\end{mytext}

\begin{mytext}
A holomorphic germ $f:(X, p) \rightarrow (Y,q)$ naturally induces an action $f^*$ on $R=\nC[[x,y]]$, by composition: $\phi \mapsto f^*\phi = \phi \circ f$.
The natural way to define an action on (centered) valuations seems to be the dual action $f_*\nu=\nu \circ f^*$; explicitly we have $f_*\nu(\phi)=\nu(\phi \circ f)$.
This definition works for Krull valuations, but not for valuations: if $\nu \in \mc{V}$, then clearly $f_*\nu$ is a valuation, but it might not be proper.
More precisely, $f_*\nu$ is not centered if and only if $\nu=\nu_C$ is a curve valuation, with $C=\{\phi=0\}$ an irreducible curve contracted to $q$ by $f$ (that is to say if $f^*\mg{m} \subseteq \langle \phi \rangle$). In this case $C$ has to be a critical curve, and $f_* \nu$ is not proper.
\end{mytext}

\begin{mydef}
Let $f:(X, p) \rightarrow (Y,q)$ be a (dominant) holomorphic germ. We call \mydefal{contracted critical curve valuations}{contracted critical curve valuations} for $f$ the valuations $\nu_C$ with $C$ a critical curve contracted to $q$ by $f$. We denote by $\mg{C}_f$ the set of all contracted critical curve valuations for $f$. 
\end{mydef}

\begin{myoss}
$\mg{C}_f$ has a finite number of elements, all ends for the valuative tree.
\end{myoss}

\begin{mytext}
So if $\nu \in \mc{V}\setminus \mg{C}_f$, then $f_*\nu$ is a centered valuation, but not normalized generally. The norm will be $f_*\nu(\mg{m})=\nu(f^*\mg{m})$: we can renormalize this valuation and obtain an action $f_\bullet: \mc{V}\setminus \mg{C}_f \rightarrow \mc{V}$.
\end{mytext}

\begin{mydef}
Let $f:(X, p) \rightarrow (Y,q)$ be a (dominant) holomorphic germ. For every valuation $\nu \in \mc{V}$ we define $c(f, \nu):= \nu(f^*\mg{m})$ the \mydefal{attraction rate of $f$ along $\nu$}{attraction rate}; if $\nu=\nu_{\mg{m}}$ is the multiplicity valuation, then we simply write $c(f):=c(f, \nu_{\mg{m}})$ the \mydefal{attraction rate}{attraction rate} of $f$.
For every valuation $\nu \in \mc{V} \setminus \mg{C}_f$ we define $f_\bullet \nu := f_*\nu / c(f,\nu) \in \mc{V}$.
If $f:(X, p) \rightarrow (X,p)$, we will also define $c_\infty(f):= \lim_{n \to \infty} \sqrt[n]{c(f^n)}$ the \mydefal{asymptotic attraction rate}{asymptotic attraction rate} of $f$.
\end{mydef}

\begin{mytext}
Up to fix coordinates in $p$ and $q$, we can consider a germ $f:(X,p) \rightarrow (Y,q)$ as a germ $f:(\nC^2,0) \rightarrow (\nC^2,0)$: from now on we will state results in the latter case, but they can be easily extended to the general case.
\end{mytext}

\begin{mytext}
In order to have an action on $\mc{V}$, we should extend $f_\bullet$ to contracted critical curve valuations.
\end{mytext}

\begin{myprop}[{\cite[Proposition 2.7]{favre-jonsson:eigenval}}]
Suppose $C$ is an irreducible curve germ such that $f(C)=\{0\}$ (i.e. $\nu_C \in \mg{C}_f$). Then $c(f, \nu_C)=\infty$. Further, the limit of $f_\bullet \nu$ as $\nu$ increases to $\nu_C$ exists, and it is a divisorial valuation that we denote by $f_\bullet \nu_C$.
It can be interpreted geometrically as follows. There exist modifications $\pi:X \rightarrow (\nC^2, 0)$ and $\pi^\prime:X^\prime \rightarrow (\nC^2, 0)$, such that $f$ lifts to a holomorphic map $\hat{f}:X \rightarrow X^\prime$ sending $C$ to a curve germ included in an exceptional component $E^\prime \in \Gamma_{\pi^\prime}^*$, for which $f_\bullet \nu_C = \nu_{E^\prime}$.
\end{myprop}

\begin{mydef}
Let $f:(\nC^2,0) \rightarrow (\nC^2,0)$ be a (dominant) holomorphic germ. For every $\nu \in \mc{V}$ we shall denote by $d(f_\bullet)_\nu:T_\nu\mc{V} \rightarrow T_{f_\bullet \nu} \mc{V}$ the \mydefal{tangent map}{tree map!tangent map} induced by $f$ at $\nu$.
We will often omit the point $\nu$ where we are considering the tangent map, and write $d(f_\bullet)_\nu = df_\bullet$.
\end{mydef}

\begin{mytext}
For other properties of the action $f_\bullet$, we refer to \cite[Sections 2 and 3]{favre-jonsson:eigenval}.
\end{mytext}

\mysubsubsect{Eigenvaluations and Basins of Attraction}

\begin{mytext}
Thanks to regularity properties of $f_\bullet$ (see \cite[Theorem 3.1]{favre-jonsson:eigenval}) and a fixed point theorem for regular tree maps (see \cite[Theorem 4.5]{favre-jonsson:eigenval}) we obtain eigenvaluations.
\end{mytext}

\begin{mythm}[{\cite[Theorem 4.2]{favre-jonsson:eigenval}}]\label{thm:eigenvaluation}
Let $f:(\nC^2, 0) \rightarrow (\nC^2,0)$ be a (dominant) holomorphic germ. Then there exists a valuation $\nu_\star \in \mc{V}$ such that $f_\bullet \nu_\star=\nu_\star$, and $c(f, \nu_\star)=c_\infty(f)=:c_\infty$.
Moreover $\nu_\star$ cannot be a contracted critical curve valuation, and neither a non-analytic curve valuation if $c_\infty > 1$.
If $\nu_\star$ is an end, then there exists $\nu_0 < \nu_\star$ (arbitrarily close to $\nu_\star$), such that $c(f, \nu_0)=c_\infty$, $f_\bullet$ preserves the order on $\{\nu \geq \nu_0\}$ and $f_\bullet \nu > \nu$ for every $\nu \in [\nu_0, \nu_\star)$.
Finally, we can find $0 < \delta \leq 1$ such that $\delta c_\infty^n \leq c(f^n) \leq c_\infty^n$ for every $n \geq 1$.
\end{mythm}

\begin{mydef}
Let $f:(\nC^2, 0) \rightarrow (\nC^2,0)$ be a (dominant) holomorphic germ.
A valuation $\nu_\star \in \mc{V}$ is called \mydefin{fixed valuation} for $f$ if $f_\bullet \nu_\star=\nu_\star$. It is called \mydefal{eigenvaluation}{eigenvaluation} for $f$ if it is a quasimonomial fixed valuation, or a fixed valuation which is a strongly attracting end (see \cite[Section 4]{favre-jonsson:eigenval}).
\end{mydef}

\begin{myoss}
In the rest of this paper, we will always consider quasimonomial eigenvaluations whenever possible. Therefore when we will say that an eigenvaluation $\nu_\star$ is an end, we implicitly state that quasimonomial eigenvaluations do not exist.
\end{myoss}

\begin{mycor}\label{cor:avoidcases}
Let $f$ be a (dominant) holomorphic germ, and let $\nu_\star$ be an eigenvaluation for $f$.
Then:
\begin{enumerate}[(i)]
\item if $c_\infty(f) > 1$ then $\nu_\star$ cannot be a non-analytic curve valuation;
\item if $c_\infty(f) = 1$ then $\nu_\star$ cannot be a quasimonomial valuation.
\end{enumerate}
\end{mycor}

\begin{mydim}
The first assertion has been already stated in \ref{thm:eigenvaluation}.

Let us suppose $c_\infty(f) = 1$. Then applying \cite[Lemma 7.7]{favre-jonsson:eigenval} to the eigenvaluation (and recalling that $c_\infty(f)=c(f, \nu_\star)$ by Theorem \ref{thm:eigenvaluation}) we obtain
$$
A(\nu_\star)=A(\nu_\star) + \nu_\star(Jf) \mbox{,}
$$
that can be satisfied only if $A(\nu_\star)=\infty$. It follows that $\nu_\star$ cannot be a quasimonomial valuation.
\end{mydim}

\begin{myprop}[{\cite[Proposition 5.2]{favre-jonsson:eigenval}}]\label{prop:basinofattraction}
Let $f$ be a (dominant) holomorphic germ, and let $\nu_\star$ be an eigenvaluation for $f$.
\begin{enumerate}[(i)]
\item If $\nu_\star$ is an end for $\mc{V}$, then for any $\nu_0 \in \mc{V}$ with $\nu_0 \leq \nu_\star$, and $\nu_0$ sufficiently close to $\nu_\star$, $f_\bullet$ maps the segment $I=[\nu_0, \nu_\star]$ strictly into itself and is order-preserving there.
Moreover, if we set $U=U(\vect{v})$, where $\vect{v}$ is the tangent vector at $\nu_0$ represented by $\nu_\star$, then $f_\bullet$ also maps the open set $U$ strictly into itself and $f_\bullet^n \rightarrow \nu_\star$ as $n \rightarrow \infty$ in $U$.
\item If $\nu_\star$ is divisorial, then there exists a tangent vector $\vect{w}$ at $\nu_\star$ such that for any $\nu_0 \in \mc{V}$ representing $\vect{w}$ and sufficiently close to $\nu_\star$, $f_\bullet$ maps the segment $I=[\nu_\star, \nu_0]$ into itself and is order-preserving there.
Moreover, if we set $U=U(\vect{v}) \cap U(\vect{w})$, where $\vect{v}$ is the tangent vector at $\nu_0$ represented by $\nu_\star$, then $f_\bullet(I) \subrelcpt I$, $f_\bullet(U) \subrelcpt U$, and $f_\bullet^n \rightarrow \nu_\star$ as $n \rightarrow \infty$ on $U$.
\item If $\nu_\star$ is irrational, then there exist $\nu_1, \nu_2 \in \mc{V}$, arbitrarily close to $\nu_\star$, with $\nu_1 < \nu_\star < \nu_2$ such that $f_\bullet$ maps the segment $I=[\nu_1, \nu_2]$ into itself. Let $\vect{v_i}$ be the tangent vector at $\nu_i$ represented by $\nu_\star$ (for $i=1, 2$), and set $U=U(\vect{v_1}) \cap U(\vect{v_2})$. Then $f_\bullet(U) \subseteq U$. Further, either $f_\bullet|_I^2=\on{id}_I$ or $f_\bullet^n \rightarrow \nu_\star$ as $n \rightarrow \infty$ on $U$.
\end{enumerate}
\end{myprop}

\mysect{Rigidification}

\mysubsect{General result}

\begin{mytext}
In this section we shall prove our main theorem (Theorem \ref{thm:rigidification}). Here there are five cases instead of the four of \cite[Theorem 5.1]{favre-jonsson:eigenval}: the new case is when we have a non-analytic curve eigenvaluation, and it arises only when we deal with $f$ having a non-nilpotent differential.
For the other cases, we refer directly to \cite[Theorem 5.1]{favre-jonsson:eigenval}.
\end{mytext}

\begin{mydim}[Proof of Theorem \ref{thm:rigidification}.]
Let $\nu_\star$ be an eigenvaluation for $f$ (that exists thanks to Theorem \ref{thm:eigenvaluation}), and suppose that it is a non-analytic curve valuation $\nu_C$.

Pick $\nu_0$ as in Proposition \ref{prop:basinofattraction}. By increasing $\nu_0$, we can suppose $\nu_0$ divisorial. Let $\pi \in \mc{B}$ be a modification such that $\nu_0=\nu_{E_0}$.
From \cite[Proposition 6.32]{favre-jonsson:eigenval} there exists a unique best approximation $\nu_E$ of $\nu_\star$ for $\pi$ (it is unique because $\nu_\star$ is an end of $\mc{V}$).
We have $\nu_0 \leq \nu_E < \nu_C$, that can be chosen arbitrarily close to $\nu_C$ (by increasing $\nu_0$).
We consider now $U=U(p)=U_{\nu_E}([\nu_\star])$.

From Proposition \ref{prop:basinofattraction} and \cite[Proposition 3.2]{favre-jonsson:eigenval}, it follows $f_\bullet U \subrelcpt U$, and the lift $\hat{f}=\pi^{-1} \circ f \circ \pi$ is holomorphic in $p$, and $\hat{f}(p)=p$. By shrinking $U(p)$, we can avoid all critical curve valuations.

It follows that $\mc{C}^\infty(\hat{f})=E$ has normal crossings.
Moreover, $E$ is contracted to $p$ by $\hat{f}$ (because $f_\bullet \nu_E > \nu_E$), $\mc{C}^\infty(\hat{f})$ is forward $\hat{f}$-invariant and $\hat{f}$ is rigid.
\end{mydim}

\begin{myoss}\label{oss:uniqueeigenval}
Studying the behavior of $\pi_\bullet$, with $\pi:(X, p) \rightarrow (\nC^2, 0)$ a modification, we see that $\pi_\bullet$ is a bijection between $\mc{V}$ and $\overline{U(p)}$.
Moreover, from the relation $\hat{f}=\pi^{-1} \circ f \circ \pi$, we see that $\pi_\bullet$ gives us a conjugation between $\hat{f}_\bullet$ and $f_\bullet|_{\overline{U(p)}}$.
So from the dynamics of $f_\bullet$ on $\overline{U(p)}$, we can obtain informations on the rigidification $\hat{f}$.
For example, when $f_\bullet^n \to \nu_\star$, then $\hat{f}$ will have a unique eigenvaluation $\pi_\bullet^{-1}\nu_\star$.
\end{myoss}

\mysubsect{Semi-superattracting case}

\begin{mytext}
In this section we deal with the semi-superattracting case, proving the uniqueness of the eigenvaluation in this case (see Theorem \ref{thm:fixedvaluniqueness}).
We shall write
$$
D_\lambda :=
\left(
\begin{array}{ll}
\lambda & 0 \\
0 & 0 
\end{array}
\right)
\mbox{.}
$$
\end{mytext}

\begin{mylem}\label{lem:liftDlambda}
Let $f$ be a (dominant) semi-superattracting holomorphic germ, such that $df_0=D_\lambda$ with $\lambda \neq 0$, and let $\pi:X \rightarrow (\nC^2,0)$ be the single blow-up in $0 \in \nC^2$, with $E:=\pi^{-1}(0)\cong \nP^1(\nC)$ the exceptional divisor. Set $p=[1:0] \in E$, and let $\hat{f}:(X,p) \rightarrow (X,p)$ be the lift of $f$ through $\pi$. Then $\hat{f}$ is a semi-superattracting holomorphic germ, and $d\hat{f}_p \cong D_\lambda$.
\end{mylem}
\begin{mydim}
Since $df_0=D_\lambda$, we have
\begin{equation}\label{eqn:fDlambda}
f(z,w)=\big(\lambda z + f_1(z,w), f_2(z,w)\big)\mbox{,}
\end{equation}
with $f_1, f_2 \in \mg{m}^2$.
In the chart $\pi^{-1}(\{z \neq 0\})$ we can choose $(u, t)$ coordinates in $p \in E$ such that
$$
(z,w)=\pi(u,t)=(u,ut) \mbox{.}
$$
So for the lift $\hat{f}=\pi^{-1}\circ f \circ \pi$ we have
$$
f \circ \pi (u,t) = \big(\lambda u + f_1(u,ut), f_2(u,ut)\big) \mbox{,}
$$
and then
$$
\hat{f}(u,t) = \left( \lambda u + f_1(u,ut), \frac{f_2(u,ut)}{\lambda u + f_1(u,ut)} \right) \mbox{.}
$$

We have $u^2 \mid f_1(u,ut),f_2(u,ut)$; if we set $\hat{f}=(g_1, g_2)$, we have
\begin{align*}
g_1(u,t) & = \lambda u \big(1 + O(u)\big) \\
g_2(u,t) & = \frac{u^2 O(1)}{\lambda u \big(1 + O(u)\big)} = \alpha u + O(u^2) \mbox{,}
\end{align*}
with $\alpha = \lambda^{-1}a_{2,0}$, if $f_2(z,w)=\sum_{i+j\geq 2}a_{i,j}z^i w^j$. 
It follows that
$$
d\hat{f}_p=
\left(
\begin{array}{cc}
\lambda & 0 \\
\alpha & 0
\end{array}
\right)
\cong D_\lambda \mbox{.}
$$
So $\hat{f}$ is a holomorphic germ with $d\hat{f}_p \cong D_\lambda$.
\end{mydim}

\begin{myprop}\label{prop:liftDlambda}
Let $f$ be a (dominant) semi-superattracting holomorphic germ such that $df_0=D_\lambda$ with $\lambda \neq 0$, $\nu_\star$ an eigenvaluation for $f$, and $(\pi,p,\hat{f})$ a rigidification obtained from $\nu_\star$ as in Theorem \ref{thm:rigidification}.
Then $d\hat{f}_p \cong D_\lambda$ and $\nu_\star=\nu_C$ is a (possibly formal) curve valuation, with $m(C)=1$.
\end{myprop}
\begin{mydim}
For proving this result, we have to follow the proof of \cite[Theorem 4.5]{favre-jonsson:eigenval} under the assumption $df_0 \cong D_\lambda$.
Starting from any $\nu_0$ (as in the proof of \cite[Theorem 4.5]{favre-jonsson:eigenval}), we take any end $\nu_0^\prime > f_\bullet \nu_0$, and we consider the induced tree map $F_0$ on $I_0=[\nu_0, \nu_0^\prime]$. Let $\nu_1$ be the (minimum) fixed point of $F_0$. Since $f_\bullet$ has no quasimonomial eigenvaluations (see Corollary \ref{cor:avoidcases}), then $\nu_1 \geq f_\bullet \nu_0$. Up to choosing $\nu_0^\prime$ such that $\nu_0^\prime \not \in d(f_\bullet)_{\nu_0} ([f_\bullet \nu_0])$, we can suppose that $\nu_1 = f_\bullet \nu_0$.

Let us apply this argument for $\nu_0=\nu_{\mg{m}}$. If $f$ is as in \eqref{eqn:fDlambda}, then $\nu_1=f_\bullet \nu_0$ is a divisorial valuation associated to an exceptional component $E_1$ obtained from the exceptional component $E_0$ of a single blow-up of $0 \in \nC^2$ only by blowing-up free points (i.e., the generic multiplicity $b(\nu_{E_1})$ of $\nu_{E_1}$ is equal to $1$): as a matter of fact, $f_*\nu_0(x)=1$ while $f_*\nu_0(\phi) \in \nN$ for every $\phi \in R$.

Applying this argument recursively (as in the proof of \cite[Theorem 4.5]{favre-jonsson:eigenval}), we get the assertion on the type of eigenvaluation.

For the result on $d\hat{f}_p$, we only have to observe that, on the proofs of Theorem \ref{thm:rigidification} and \cite[Theorem 5.1]{favre-jonsson:eigenval} in the case of an analytic curve eigenvaluation, up to shrink the basin of attraction, we can choose the infinitely near point $p$ such that $\nu_{E_p}$ has generic multiplicity $b(\nu_{E_p})$ equal to $1$, where $E_p$ denotes the exceptional component obtained blowing-up $p$.
Then, the modification $\pi$ on the rigidification is the composition of blow-ups of free points, and then we can apply (recursively) Lemma \ref{lem:liftDlambda} and obtain the thesis. 
\end{mydim}

\begin{mylem}\label{lem:preparationDlambda}
Let $f$ be a (dominant) semi-superattracting holomorphic germ such that $df_0=D_\lambda$ with $\lambda \neq 0$. Then, up to a (possibly formal) change of coordinates, we can suppose that
$$
f(z,w)=\Big(\lambda z \big(1+ f_1(z,w)\big), w f_2(z,w)\Big)\mbox,
$$
with $f_1, f_2 \in \mg{m}$.
\end{mylem}
\begin{mydim}
First of all, we can suppose that
$$
f(z,w)=\big(\lambda z + g_1(z,w), g_2(z,w)\big)\mbox,
$$
with $g_1,g_2 \in \mg{m}^2$.

Thanks to Proposition \ref{prop:liftDlambda}, we know that there is an eigenvaluation $\nu_\star=\nu_C$ with $C=\{\phi=0\}$ a (possibly formal) curve, with $\phi(z,w)=w-\theta(z)$ for a suitable $\theta$.
Up to the (possibly formal) change of coordinates $(z,w) \mapsto (z,w-\theta(z))$, we can suppose that $\phi=w$, and in particular, since $C$ is fixed by $f$, that $w | g_2$.
Then we have
$$
f(z,w)=\Big(\lambda z\big(1 + f_1(z,w)\big) + h(w), w f_2(z,w)\Big)\mbox,
$$
with $f_1, f_2 \in \mg{m}$ and $h \in \mg{m}^2$.
We shall denote $g_2(z,w)=w f_2(z,w)$.

Now we only have to show that up to a (possibly formal) change of coordinates, we can suppose $h \equiv 0$.
We consider a change of coordinates of the form $\Phi(z,w)=(z+\eta(w), w)$, with $\eta \in \mg{m}^2$.
In this case we have $\Phi^{-1}(z,w)=(z-\eta(w), w)$.
So we have
\begin{align}
\Phi^{-1}\circ f \circ \Phi(z,w)& \label{eqn:lemconj1}\\
=\Big(\lambda \big(z+\eta(w)\big)&\big(1+f_1 \circ \Phi(z,w)\big) +h(w) - \eta \circ g_2 \circ \Phi (z,w),w f_2 \circ \Phi(z,w)\Big)\mbox.\nonumber
\end{align}
We notice that the second coordinate of \eqref{eqn:lemconj1} is always divisible by $w$; we only have to show that there exists a suitable $\eta$ such that the first coordinate of \eqref{eqn:lemconj1}, valuated on $(0,w)$, is equal to $0$.
Hence we have to solve
\begin{equation}\label{eqn:lemconjtosolve}
\lambda \eta(w) \Big(1+f_1\big(\eta(w),w\big)\Big)+h(w)-\eta \circ g_2 \big(\eta(w),w\big)=0\mbox.
\end{equation}
If we set $\eta(w)=\sum_{n\geq 2} \eta_n w^n$, $h(w)=\sum_{n\geq 2} h_n w^n$, $1+f_1(z,w)=\sum_{i+j \geq 0} f_{i,j} z^i w^j$ and $g_2(z,w)=\sum_{i+j \geq 2} g_{i,j}z^i w^j$ (with $g_{n,0}=0$ for every $n$), then we have
\begin{equation}\label{eqn:lemconjric}
\lambda \sum_{i+j \geq 0} f_{i,j} \sum_{H \in \nN^{i+1}} \eta_{H} w^{\abs{H}+j}+ \sum_{n\geq 2} h_n w^n= \sum_{k} \eta_k \left(\sum_{i+j\geq 2}g_{i,j} \eta(w)^i w^j\right)^k \mbox.
\end{equation}
Comparing the coefficients of $w^n$ in both members, we get
$$
\lambda \eta_n + \mbox{l.o.t.} = \mbox{l.o.t.,}
$$
where l.o.t. denotes a suitable function depending on $\eta_h$ only for $h < n$.
So thanks to \eqref{eqn:lemconjric} we have a recurrence relation for the coefficients $\eta_n$ that is a solution of \eqref{eqn:lemconjtosolve}.
\end{mydim}

\begin{mydim}[Proof of Theorem \ref{thm:fixedvaluniqueness}.]
Thanks to Lemma \ref{lem:preparationDlambda}, we can suppose (up to formal conjugacy) that
$$
f(z,w)=\Big(\lambda z \big(1+g_1(z,w)\big), w g_2(z,w)\Big)\mbox,
$$
where $g_1, g_2 \in \mg{m}$. We shall denote $f_2(z,w)=w g_2(z,w)$.

It follows that the eigenvaluation $\nu_\star$ given by Proposition \ref{prop:liftDlambda} is $\nu_\star=\nu_w$, while $\nu_z$ is either fixed by $f_\bullet$ or a contracted critical curve valuation.

We only have to show that there are no other fixed valuations.

First of all, we want to notice what happens during the process used in the proof of Proposition \ref{prop:liftDlambda} to tangent vectors at the valuation $\nu_0=\nu_\mg{m}$.
Let us consider the family of valuations $\nu_{\theta, t}$, where $\theta \in \nP^1(\nC)$ and $t \in [1, \infty]$, described as follows:
if we denote $\phi_\theta=w-\theta z$ when $\theta \in \nC$, and $\psi_\infty=z$, then $\nu_{\theta,t}$ is the valuation of skewness $\alpha(\nu_{\theta,t})=t$ in the segment $[\nu_{\mg{m}}, \nu_{\phi_\theta}]$, i.e., the monomial valuation defined by $\nu_{\theta,t}(\phi_\theta)=t$ and $\nu_{\theta,t}(z)=1$ if $\theta \in \nC$, and $\nu_{\infty,t}(z)=t$, $\nu_{\infty,t}(w)=1$.

Then we have that $\nu_1=f_\bullet(\nu_{\mg{m}})=\nu_{0,m(f_2)}$, where $m$ denotes the multiplicity function, while $f_\bullet(\nu_{\theta,t}) \geq \nu_1$ for every $\theta \in \nC$ and $t$, since $f_\bullet(\nu_{\theta,t})(z)=1=\nu_1(z)$, $f_\bullet(\nu_{\theta,t})(w) \geq m(f_2) = \nu_1(w)$ and $\nu_1$ is the minimum valuation that assumes those values on $z$ and $w$.

We shall denote by $\vect{v_\theta}$ the tangent vector in $\nu_\mg{m}$ represented by $\nu_{\theta, \infty}$, and by $\vect{u_\infty}$ the tangent vector in $\nu_1$ represented by $\nu_{\mg{m}}$; then it follows from what we have seen that $df_\bullet(\vect{v_\theta})\neq \vect{u_\infty}$ for every $\theta \neq \infty$, and hence there are no fixed valuations in $\overline{U_{\nu_\mg{m}}(\vect{v_\theta})}$ for every $\theta \neq 0, \infty$.

Moreover, applying this argument recursively as in the proof of Proposition \ref{prop:liftDlambda}, we obtain that there are no other fixed valuations in $\overline{U_{\nu_\mg{m}}(\vect{v_0})}$, except for the eigenvaluation $\nu_w$.

It remains to check for valuations in $\overline{U_{\nu_\mg{m}}(\vect{v_\infty})}$.
For this purpose, let us consider $f_\bullet(\nu_{\infty, t})$.
For simplicity, we shall denote $\nu_{\infty,t}=\nu_{0,1/t}$ for every $t \in [0,1]$
From direct computation, we have that $f_*(\nu_{\infty, t})(z)=t$, while 
$$
f_*(\nu_{\infty, t})(w)=\bigwedge_j (a_j t +b_j)\mbox,
$$
for suitable $a_j \in \nN^*, b_j \in \nN$.
It follows that $f_\bullet(\nu_{\infty, t})=\nu_{\infty, g(t)}$ for a suitable map $g(t)$, such that $g(t) < t$, and that $d(f_\bullet)_{\nu_{\infty, t}}([\nu_w])=[\nu_w]$ (where the latter tangent vector belongs to the proper tangent space).
Letting $t$ go to $\infty$, we obtain that the only fixed valuation in $U_{\nu_\mg{m}}(\vect{v_\infty})$ is $\nu_z$, and we are done. 
\end{mydim}

\begin{myoss}\label{oss:fixedcurves}
Theorem \ref{thm:fixedvaluniqueness} shows that every semi-superattracting germ $f$ has two (formal) invariant curves: the first one $C$ associated to the eigenvaluation, and hence to the eigenvalue $\lambda$ of $df_0$; the second one $D$ associated to the fixed or contracted critical curve valuation, and hence to the eigenvalue $0$ of $df_0$.
If $f$ is of type $(0, \nC \setminus \overline{\nD})$, then both these curves are actually holomorphic, thanks to the Stable/Unstable Manifold Theorem (see \cite[Theorem 3.1.2 and Theorem 3.1.3]{abate:hypdynsys}).
In the general case of $f$ of type $(0, \nC^*)$, one can at least recover the manifold associated to the eigenvalue $0$ of $df_0$, using generalisations of the Stable Manifold Theorem, such as the Hadamard-Perron Theorem (see \cite[Theorem 3.1.4]{abate:hypdynsys}).
In particular the curve $D$ is always holomorphic.
However $C$ is not always holomorphic in general (see for example Proposition \ref{prop:rigidattr}).
\end{myoss}

\mysect{Rigid Germs}

\begin{mytext}
In this section we will introduce the classification of attracting rigid germs in $(\nC^2, 0)$ up to holomorphic and formal conjugacy (for proofs, see \cite{favre:rigidgerms}), and the classification of rigid germs of type $(0, \nC \setminus \nD)$ in $(\nC^2,0)$ up to formal conjugacy.
\end{mytext}

\begin{mytext}
For stating them, we shall need $3$ invariants.
\begin{itemize}
\item \textbf{The generalized critical set}: if $f:(\nC^2, 0) \rightarrow (\nC^2,0)$ is a rigid germ then $C=\mc{C}^\infty(f)$ is a curve with normal crossings at the origin, and it can have $0$, $1$ or $2$ irreducible components, that is to say that $\mc{C}^\infty(f)$ can be empty (if and only if $f$ is a local biholomorphism in $0$), an irreducible curve, or a reducible curve (with only $2$ irreducible components); we will call $f$ \mydefal{regular}{rigid germ!regular}, \mydefal{irreducible}{rigid germ!irreducible} or \mydefal{reducible}{rigid germ!reducible} respectively.
\item \textbf{The trace}: if $f$ is not regular, we have $2$ cases: either $\on{tr}df_0 \neq 0$, and $df_0$ has a zero eigenvalue and a non-zero eigenvalue, or $\on{tr}df_0 = 0$, and $df_0$ is nilpotent.
\item \textbf{The action on $\pi_1(\Delta^2 \setminus \mc{C}^\infty(f))$}: as $\mc{C}^\infty(f)$ is backward invariant, $f$ induces a map from $U=\Delta^2 \setminus \mc{C}^\infty(f)$ (here $\Delta^2$ denotes a sufficiently small polydisc) to itself, and so an action $f_*$ on the first fundamental group of $U$. When $f$ is irreducible, then $\pi_1(U) \cong \nZ$, and $f_*$ is completely described by $f_*(1) \in \nN^*$ ($f$ preserves orientation); when $f$ is reducible, then $\pi_1(U) \cong \nZ \oplus \nZ$, and $f_*$ is described by a $2 \times 2$ matrix with integer entries (in $\nN$).
\end{itemize}
\end{mytext}

\begin{mydef}\label{def:rigidclasses}
Let $f:(\nC^2, 0) \rightarrow (\nC^2,0)$ a rigid germ. Then $f$ belongs to:
\begin{enumerate}[\textbf{Class} $1$]
\item if $f$ is regular;
\item if $f$ is irreducible, $\on{tr} df_0 \neq 0$ and $f_*(1)=1$;
\item if $f$ is irreducible, $\on{tr} df_0 \neq 0$ and $f_*(1)\geq 2$;
\item if $f$ is irreducible, $\on{tr} df_0 = 0$ (this implies $f_*(1)\geq 2$);
\item if $f$ is reducible, $\on{tr} df_0 \neq 0$ (this implies $\det f_* \neq 0$);
\item if $f$ is reducible, $\on{tr} df_0 = 0$ and $\det f_* \neq 0$;
\item if $f$ is reducible, $\on{tr} df_0 = 0$ and $\det f_* = 0$.
\end{enumerate}
\end{mydef}

\begin{tabular}[htb]{|c|c|c|c|}
\hline
Class & $\mc{C}^\infty(f)$ & $\on{tr}df_0$ & $\det f_*$ \\
\hline
$1$ & $0$ (empty) & & \\
$2$ & $1$ (irreducible) & $\neq 0$ & $=1$ \\
$3$ & & & $\geq 2$ \\
$4$ & & $=0$ & ($\geq 2$) \\
$5$ & $2$ (reducible) & $\neq 0$ & ($\neq 0$) \\
$6$ & & $=0$ & $\neq 0$ \\
$7$ & & & $= 0$ \\
\hline
\end{tabular}

\begin{myoss}\label{oss:normfrigid}
If $f$ is irreducible, up to a change of coordinates we can assume $\mc{C}^\infty(f)=\{z=0\}$. Then just using that $\{z=0\}$ is backward invariant, we can write $f$ in the form
$$
f(z,w)=\Big(\alpha z^p \big(1+ \phi(z,w)\big), f_2(z,w)\Big) \mbox{,}
$$
with $\phi, f_2 \in \mg{m}$.
It can be easily seen that $f_*=p \geq 1$.

Analogously, if $f$ is reducible, up to a change of coordinates we can assume $\mc{C}^\infty(f)=\{zw=0\}$. Then just using that $\{zw=0\}$ is backward invariant we can write $f$ in the form
$$
f(z,w)=\Big(\lambda_1 z^a w^b \big(1+ \phi_1(z,w)\big), \lambda_2 z^c w^d \big(1+ \phi_2(z,w)\big)\Big) \mbox{,}
$$
with $\phi_1, \phi_2 \in \mg{m}$.
In this case $f_*$ is represented by the $2 \times 2$ matrix
\begin{equation}\label{eqn:Mf}
M(f):=
\left(
\begin{array}{ll}
a & b \\
c & d 
\end{array}
\right)
\mbox{.}
\end{equation}
\end{myoss}

\mysubsect{Attracting rigid germs}

\begin{mytext}
The classification up to holomorphic conjugacy of attracting rigid germs in $\nC^2$ is given in \cite[Ch.1]{favre:rigidgerms}; the only remark needed is the following.
\end{mytext}

\begin{myoss}
During the proof of \cite[Step 1 on page 491, and First case on page 498]{favre:rigidgerms}, the author starts from a germ of the form
$$
f(z,w)=\Big(\alpha z^p \big(1+ g(z,w)\big), f_2(z,w)\Big) \mbox{,}
$$
with $\phi, f_2 \in \mg{m}$, and uses K{\oe}nigs and B{\"o}ttcher Theorems (see \cite[Theorem 3.1 and Theorem 3.2]{favre:rigidgerms} respectively, \cite{koenigs:1884} and \cite{bottcher:1904} for original papers, and \cite[Theorem 8.2 and Theorem 9.1]{milnor:dyn1cplxvar} for a modern exposition of proofs) to assume, up to holomorphic conjugacy, that $g \equiv 0$ (and $\alpha=1$ if $p \geq 2$).
That argument does not work.
Let us denote by $\Phi(z,w)=(\phi_w(z),w)$ the conjugation given by those theorems, and $\tilde{f}=\Phi \circ f \circ \Phi^{-1}$. We shall also denote $f(z,w)=(f^{(1)}_w(z),f^{(2)}_w(z))$, and analogously for $\tilde{f}$.
By hypothesis $\phi_w(z)$ is such that $\phi_w \circ f^{(1)}_w \circ \phi_w^{-1}(z)=\alpha z^p$ (with $\alpha=1$ if $p\geq 2$).

Then we have that
$$
\tilde{f}^{(1)}_w(z)= \phi_{f^{(2)}_w(\phi_w^{-1}(z))} \circ f^{(1)}_w \circ \phi_w^{-1}(z)\mbox{,}
$$
and hence it does not coincide with $\alpha z^p$.

We also note that if $\abs{\alpha} > 1$, there still is a K{\oe}nigs Theorem, but the result is false (see Counterexample \ref{ces:nonhol1}).

Nevertheless, one can obtain this result in the attracting case in the following way.

We want to solve the conjugacy relation
\begin{equation}\label{eqn:conjrelation}
\Phi \circ f = e \circ \Phi\mbox{,}
\end{equation}
where $e$ is a germ of the form
$$
e(z,w)=\big(\alpha z^p, e_2(z,w)\big) \mbox{,}
$$
with $e_2 \in \mg{m}$.

We look for a solution of the form
$$
\Phi(z,w)=\Big(z\big(1+\phi(z,w)\big), w\Big)\mbox,
$$
with $\phi \in \mg{m}$.

Then from the conjugacy relation \eqref{eqn:conjrelation} (comparing the first coordinate) we get
$$
(1+g)(1+\phi \circ f)=(1+\phi)^p\mbox.
$$
Then we can consider
\begin{equation}\label{eqn:defphip1}
1+\phi=\prod_{k=0}^\infty \left(1+g \circ f^k\right)^{1/p^{k+1}}\mbox,
\end{equation}
that would work if that product converges.
But since $f$ is attracting, there exists $0 < \varepsilon < 1$ such that $\norm{f(z,w)} \leq \varepsilon \norm{(z,w)}$, while since $g \in \mg{m}$, there exists $M > 0$ such that $\abs{g(z,w)} \leq M \norm{(z,w)}$. It follows that
$$
\sum_{k=0}^\infty p^{-(k+1)} \abs{g \circ f^k(z,w)} \leq \sum_{k=0}^\infty \frac{M}{p} \left(\frac{\varepsilon}{p}\right)^k = \frac{M}{p-\varepsilon} < \infty\mbox,
$$
and hence \eqref{eqn:defphip1} defines an holomorphic germ $\phi$, and hence a holomorphic map $\Phi$ that satisfies the conjugacy relation \eqref{eqn:conjrelation} in the first coordinate.

To choose $e_2$ such that \eqref{eqn:conjrelation} holds also for the second coordinate, we have to solve
$$
f_2 = e_2 \circ \Phi\mbox,
$$
but since $\Phi$ is a holomorphic invertible map, we can just define $e_2=f_2 \circ \Phi$, and we are done.

Notice that this approach would not work for rigid germs of type $(0,\nC \setminus \nD)$, not even formally.
\end{myoss}

\mysubsect{Rigid germs of type $(0, \nC \setminus \nD)$}

\begin{mytext}
Here we are going to study formal normal forms for rigid germs of type $(0, \nC \setminus \nD)$.
As notation, if $f(z,w)=\sum_{i,j}f_{i,j}z^i w^j$ is a formal power series, and $I=(i_1, \ldots, i_k)$ and $J=(j_1, \ldots, j_k)$ are two multi-indices, then we shall denote by $f_{I,J}$ the product
$$
f_{I,J}=\prod_{l=1}^k f_{i_l, j_l}\mbox.
$$
Moreover, when writing the dummy variables of a sum, we shall write the dimension of a multi-index after the multi-index itself. For example, $I(n)$ shall denote a multi-index $I \in \nN^n$.
We shall group together the multi-indices with the same dimension, separating these groups by a semi-colon, and we shall omit the dimension when it is equal to $1$. For example,
$$
\sum_{n,m;I,J(n);K(m)}
$$
shall denote a sum over $n,m \in \nN$, $I,J \in \nN^n$ and $K \in \nN^m$.
As a convention, a multi-index of dimension $0$ is an empty multi-index.
\end{mytext}

\begin{mytext}
First of all, we need to recall the formal classification of (invertible) germs in one complex variable (for the proof, and standard theory of dynamics in one complex variable, we refer to \cite{milnor:dyn1cplxvar}).
\end{mytext}

\begin{myprop}[Formal classification in $(\nC,0)$]\label{prop:forclass1d}
Let $f:(\nC,0) \rightarrow (\nC,0)$ be a holomorphic germ, and denote by $\lambda=f^\prime(0)$ the multiplier.
Then
\begin{enumerate}[(i)]
\item if $\lambda=0$, then $f$ is formally conjugated to $z \mapsto z^p$ for a suitable $p \geq 2$;
\item if $\lambda \neq 0$, and $\lambda^r \neq 1$ for any $r \in \nN^*$, then $f$ is formally conjugated to $z \mapsto \lambda z$;
\item if $\lambda^r = 1$, then there exist (unique) $s \in r\nN^*$ and $\beta \in \nC$ such that $f$ is formally conjugated to $z \mapsto z(1 + z^s + \beta z^{2s})$.
\end{enumerate}
\end{myprop}

\begin{myoss}\label{oss:firstaction}
If $f:(\nC^2,0)\rightarrow (\nC^2,0)$ is a semi-superattracting holomorphic germ, recalling Remark \ref{oss:fixedcurves}, we have two invariant curves, $C$ and $D$, with transverse intersection and multiplicity equal to $1$, that play the role of the Unstable-Stable manifold. In particular, the formal conjugacy classes of $f|_C$ and $f|_D$ are formal invariants.

Moreover, up to formal conjugacy, we can suppose that $C=\{w=0\}$ and $D=\{z=0\}$. Let us set $f=(f_1,f_2)$; then, up to a formal change of coordinates, we can suppose that $f_1(z,0)$ is equal to one of the formal normal forms given by Proposition \ref{prop:forclass1d}.

Indeed, if $\phi \in \nC[[z]]$ is the formal conjugation between $f_1(z,0)$ and its formal conjugacy class $h(z)$, the formal map $\Phi(z,w)=(\phi(z),w)$ is a conjugation between $f$ and a map $g$, with $g_1(\cdot,0)=h(\cdot)$.

We shall refer at the normal form $h$ of a germ $f$ as the \mydefin{first (formal) action} of $f$.
\end{myoss}

\begin{mylem}\label{lem:firstaction}
Let $f:(\nC^2,0)\rightarrow (\nC^2,0)$ be a semi-superattracting holomorphic germ. Then, up to formal conjugacy, we can suppose that
$$
f(z,w)=\big(h(z), g(z,w)\big) \mbox,
$$
with $h$ the first action of $f$, and $g \in \mg{m}^2$.
\end{mylem}
\begin{mydim}
We can suppose that $f$ is of the form
$$
f(z,w)=\Big(\lambda z \big(1+f_1(z,w)\big), g_2(z,w)\Big)\mbox,
$$
with $f_1 \in \mg{m}$ and $w | g_2 \in \mg{m}^2$.

We want to find a conjugation map of the form $\Phi(z,w)=(z(1+\phi(z,w)), w)$ that conjugates $f$ with
$$
e(z,w)=\Big(\lambda z \big(1+e_1(z)\big),e_2(z,w)\Big)\mbox,
$$
with $w | e_2 \in \mg{m}^2$, and $\lambda z (1+e_1(z))=h(z)$. 

Let us set $1+f_1(z,w)=\sum_{i+j\geq 0} f_{i,j} z^i w^j$, $g_2(z,w)=\sum_{i+j \geq 2} g_{i,j} z^i w^j$, $1+\phi(z,w)=\sum_{i+j \geq 0} \phi_{i,j} z^i w^j$ and $1+e_1(z)=\sum_{i\geq 0} e_{i} z^i$.

Then for the first coordinate of the conjugacy equation $\Phi \circ f = e \circ \Phi$ we have
\begin{align}
\sum_{i+j \geq 0} \phi_{i,j} \lambda^{i+1} z^{i+1} \sum_{I,J \in \nN^{i+1}}& f_{I,J} z^{\abs{I}} w^{\abs{J}} \sum_{H,K \in \nN^j} g_{H,K} z^{\abs{H}} w^{\abs{K}} \label{eqn:conj1lhs}\\ 
&|\,| \label{eqn:conj1} \\
\lambda \sum_{h} e_h z^{h+1} \sum_{N,M \in \nN^{h+1}} & \phi_{N,M} z^{\abs{N}} w^{\abs{M}}\label{eqn:conj1rhs}\mbox.
\end{align}
If we denote by $\lhs_{n,m}$ and by $\rhs_{n,m}$ the coefficients of $z^n w^m$ respectively of \eqref{eqn:conj1lhs} and \eqref{eqn:conj1rhs}, we have 
\begin{align*}
\lhs_{n,m} =& \sumaa{\substack{i,j;I,J (i+1); H,K (j)\\ i+1+\abs{I}+\abs{H}=n\\ \abs{J}+\abs{K}=m}}{0.8cm} \phi_{i,j} \lambda^{i+1} f_{I,J} g_{H,K}\mbox; &
\rhs_{n,m} =& \sumaa{\substack{h;N,M (h+1)\\ h+1+\abs{N}=n\\ \abs{M}=m}}{0.4cm} \lambda e_h \phi_{N,M}\mbox.
\end{align*}
If we denote by lower order terms all terms depending on $\phi_{i,j}$ for $(i,j)$ lower than the ones that compares in the equation (with respect to the lexicographic order), we get
$$
\delta_m^0 \phi_{n-1,0} \lambda^n f_{0,0}^n + \text{l.o.t.}= \lhs_{n,m} = \rhs_{n,m} = \lambda e_0 \phi_{n-1,m} + \text{l.o.t.}
$$ 
In particular, for $n=0$ we have $0=\lhs_{0,m} = \rhs_{0,m} = 0$ for every $m \in \nN$, while for every $m \geq 1$ we have $\lhs_{n,m}=\text{l.o.t.}$ for every $n \in \nN^*$. Since $\lambda e_0 = \lambda \neq 0$, we can use \eqref{eqn:conj1} to define recursively $\phi_{n,m}$ for every $m \geq 1$ once we have defined the base step for $m=0$.

But the case $m=0$ is exactly the same as consider the formal classification of $\tilde{f}(z)=\lambda z (1+f_1(z,0))$ as a map in one complex variable. Then, again recalling Remark \ref{oss:firstaction} and putting all together, we can define a formal map $\Phi$ that solves the conjugacy relation $\Phi \circ f = e \circ \Phi$.
\end{mydim}

\begin{mydim}[Proof of Theorem \ref{thm:rigidsemiattr}.]
Thanks to Lemma \ref{lem:firstaction} and simple considerations on rigid germs (see Remark \ref{oss:normfrigid}), we can suppose that
$$
f(z,w)=\Big(h(z),z^c w^d \big(1+g(z,w)\big)\Big)
$$
for a suitable $g \in \mg{m}$, and where $h(z)=\lambda z (1+\delta(z))$ is the first action of $f$.

We want to find a conjugation $\Psi$ of the form $\Psi(z,w)=(z,w (1+\psi(z,w)))$, between $f$ and
$$
e(z,w)=\Big(h(z),z^c w^d \big(1+\varepsilon(z)\big)\Big)\mbox,
$$
for a suitable $\varepsilon$.

Let us set $\delta(z)=\sum_{i\geq 0} \delta_i z^i$, $g(z,w)=\sum_{i+j \geq 2} g_{i,j} z^i w^j$, $1+\phi(z,w)=\sum_{i+j \geq 0} \phi_{i,j} z^i w^j$ and $1+\varepsilon(z)=\sum_{i \geq 0} \varepsilon_i z^i$.

Then for the second coordinate of the conjugacy equation $\Psi \circ f = e \circ \Psi$ we have
\begin{align}
z^c w^d \sum_{i+j \geq 0} \psi_{i,j} \lambda^{i} z^{i} \sum_{L \in \nN^{i}}& \delta_{L} z^{\abs{L}} \sum_{I,J \in \nN^{j+1}} g_{I,J} z^{\abs{I}} w^{\abs{J}} \label{eqn:conj2lhs}\\
&|\,| \label{eqn:conj2} \\
z^c w^d \sum_{h} \varepsilon_{h} z^{h} \sum_{H,K \in \nN^{d}} & \psi_{H,K} z^{\abs{H}} w^{\abs{K}}\label{eqn:conj2rhs}\mbox.
\end{align}

If we denote by $\lhs_{n,m}$ and by $\rhs_{n,m}$ the coefficients of $z^{c+n} w^{d+m}$ respectively of \eqref{eqn:conj2lhs} and \eqref{eqn:conj2rhs}, we have 
\begin{align*}
\lhs_{n,m} =& \sumaa{\substack{i,j;L(i); I,J (j+1)\\ i+cj+\abs{L}+\abs{I}=n\\ dj + \abs{J}=m}}{0.6cm} \psi_{i,j} \lambda^{i} \delta_{L} g_{I,J}\mbox; &
\rhs_{n,m} =& \sumaa{\substack{h;H,K (d)\\ h+\abs{H}=n\\ \abs{K}=m}}{0.4cm} \lambda \varepsilon_{h} \psi_{H,K}\mbox;
\end{align*}
Then for $(n,m)\neq (0,0)$ we get
$$
\delta_m^0 \psi_{n,0} \lambda^n + \text{l.o.t.}= \lhs_{n,m} = \rhs_{n,m} = d \psi_{n,m} + \text{l.o.t.}
$$
Hence if $m > 0$, we can use \eqref{eqn:conj2} to define recursively $\psi_{n,m}$; for $m=0$, we can have some resonancy problems, when $\lambda^n=d$, that is exactly the condition expressed on (ii) and (iii).
In these cases, studying the dependence of $\rhs_{n,0}$ on $\varepsilon_h$, we get
$$
\rhs_{n,0} = \varepsilon_n + \text{l.o.t.}\mbox,
$$
where l.o.t. denotes here the dependence on lower order terms $\varepsilon_h$ with $h < n$.
So for each $n$ that gives us a resonance, there exists a $\varepsilon_n$ that satisfy $\lhs_{n,0}=\rhs_{n,0}$.
Putting all together, and eventually performing a conjugacy by a linear map, we obtain the thesis.
\end{mydim}

\begin{myoss}
We notice that in the statement of Theorem \ref{thm:rigidsemiattr}, $f$ belongs to Class $2$ if and only if $d=1$, to Class $3$ if and only if $c=0$, and to Class $5$ otherwise.
\end{myoss}

\begin{myoss}\label{oss:resonancebulletaction}
The composition $\alpha \circ f_\bullet$, where $\alpha$ is either skewness or thinness, is not affected by slightly changing the non-null coefficients of a germ $f$ (as far as we keep these coefficients non-null). What changes is the action of the differential $df_\bullet$ in suitable tangent spaces.

So the difference between normal forms in the resonant case of Theorem \ref{thm:rigidsemiattr} lies in the action of $df_\bullet$, that is not invariant by change of coordinates, but has a very complicated behavior.
\end{myoss}

\begin{myoss}\label{oss:convergenceseries}
Let $\phi(z,w)=\sum \phi_{n,m} z^n w^m$ be a formal power series. Then $\phi$ is holomorphic (as a germ in $0$) if and only if there is $M$ such that
$$
\abs{\phi_{n,m}} \leq M \alpha^n \beta^m\mbox.
$$
In particular, if $\phi$ is holomorphic, then $\limsup_n \sqrt[n]{\abs{\phi_{n,m}}} < \infty$ for every $m \in \nN$, and the same holds if we exchange the role of $m$ and $n$.
\end{myoss}

We shall see now that when one has rigid germs of type $(0,\nC \setminus \overline{\nD})$, one cannot generally perform either the conjugacy of Lemma \ref{lem:firstaction} or the one of Theorem \ref{thm:rigidsemiattr} in a holomorphic way (this behavior is the opposite of the $(0, \nD)$ case).

\begin{myces}\label{ces:nonhol1}
Let us show that the conjugation given by Lemma \ref{lem:firstaction} cannot be always holomorphic.
Let $f(z,w)=(\lambda z (1+w), zw)$ with $\abs{\lambda}>1$ and $e(z,w)=(\lambda z, e_2)$, and let
$$
\Phi(z,w)=\Big(z\big(1+\phi(z,w)\big), \psi(z,w)\Big)
$$
be the (formal) conjugation given by Lemma \ref{lem:firstaction}.
By direct computation, we get
$$
\phi_{n,1}=\lambda^{n(n-1)/2}\mbox.
$$
Recalling Remark \ref{oss:convergenceseries}, we have that $\phi$ is not holomorphic.
\end{myces}

\begin{myoss}
Suppose that we have a germ $f=(\lambda z, f_2)$, with $\abs{\lambda}>1$, that is formally conjugated to $(\lambda z, z^c w^d)$ (i.e., we are not in the resonance case). The proof of Theorem \ref{thm:rigidsemiattr} shows also that the conjugation with the normal forms is  unique when we ask it to be of the form $\Psi(z,w)=(z,w(1+\psi))$. 
But if we consider a general conjugation map $\Phi=(\phi_1, \phi_2)$, since we have two invariant curves $D=\{z=0\}$ and $C=\{w=0\}$, then we have $z | \phi_1$ and $w | \phi_2$, and since the first coordinate is $\lambda z$, from direct computation we also have that $\phi_1(z,w)=z$. So $\Phi$ is unique up to a linear change of coordinates.
Hence, to prove that two germ are formally but not holomorphically conjugated, we only have to show that the conjugation found during the proof of Theorem \ref{thm:rigidsemiattr} is not holomorphic.
\end{myoss}

\begin{myces}\label{ces:nonhol2}
Let us show that also the conjugation given by Theorem \ref{thm:rigidsemiattr} cannot be always holomorphic.
Let $f(z,w)=(\lambda z, zw(1+w))$ with $\abs{\lambda} > 1$ and $e(z,w)=(\lambda z, zw)$, and let
$$
\Psi(z,w)=\Big(z, w\big(1+\psi(z,w)\big)\Big)
$$
be the (formal) conjugation given by Theorem \ref{thm:rigidsemiattr}.
By direct computation, we get
$$
\psi_{n,1}=\lambda^{n(n-1)/2}\mbox,
$$
and again we have that $\psi$ is not holomorphic.
\end{myces}

\mysect{Normal Forms}

\mysubsect{Nilpotent case}

\begin{mytext}
Favre and Jonsson studied the superattracting case (see \cite[Theorem 5.1]{favre-jonsson:eigenval}); the nilpotent case is almost the same, there is in fact just one little difference between them: to explain it, we first just prove this simple Lemma.
\end{mytext}

\begin{mylem}\label{lem:liftcinfty}
Let $f$ be a (dominant) holomorphic germ, with $df_0$ non-invertible, $\nu_\star$ an eigenvaluation for $f$, and $(\pi,p,\hat{f})$ a rigidification obtained from $\nu_\star$ as in Theorem \ref{thm:rigidification}.

Assume $\nu_\star$ is not a divisorial valuation. Then $c_\infty(\hat{f})=c_\infty(f)$.
\end{mylem}
\begin{mydim}
Directly from the definition of $\hat{f}$ as lift of $f$, we have $\pi \circ \hat{f}=f \circ \pi$. Let $\mu_\star= \pi_\bullet^{-1} (\nu_\star)$ (in this case $\mu_\star$ is an eigenvaluation for $\hat{f}$). Then:
\begin{align*}
c(\pi \circ &\hat{f},\mu_\star) = c(\hat{f},\mu_\star) \cdot c(\pi, \hat{f}_\bullet \mu_\star) = c(\hat{f},\mu_\star) \cdot c(\pi, \mu_\star) \\
&|| \\
c(f \circ &\pi,\mu_\star) = c(\pi,\mu_\star) \cdot c(f, \pi_\bullet \mu_\star) = c(\pi,\mu_\star) \cdot c(f, \nu_\star)\mbox{.}
\end{align*}
From Theorem \ref{thm:eigenvaluation} we have $c(f, \nu_\star)=c_\infty(f)$ and $c(\hat{f},\mu_\star)=c_\infty(\hat{f})$; so, if $c(\pi,\mu_\star) < \infty$, we have $c_\infty(f)=c_\infty(\hat{f})$.
But $c(\pi,\mu_\star) = \infty$ if and only if $\mu_\star \in \partial U(p)$; following the proof of Theorem \ref{thm:rigidification}, this (always) happens if and only if $\nu_\star$ is a divisorial valuation.
\end{mydim}

\begin{myoss}\label{oss:superattrvsnilpot}
The unique difference between the superattracting case and the nilpotent case is that, in the nilpotent case, one has $c_\infty(f)\geq \sqrt{2}$, while in the superattracting case one has $c_\infty(f)\geq 2$.
Moreover, thanks to Lemma \ref{lem:liftcinfty}, when the eigenvaluation $\nu_\star$ is not divisorial, then for the lift $\hat{f}$ we have $c_\infty(f)=c_\infty(\hat{f})$.
So to obtain the result for the nilpotent case, we have just to ignore the hypothesis $c_\infty(\hat{f}) \geq 2$ (when $\nu_\star$ is not divisorial).
\end{myoss}

\mysubsubsect{Germs of type $(0, \nD^*)$}

\begin{myprop}\label{prop:rigidattr}
Let $f$ be a (dominant) holomorphic germ of type $(0, \nD^*)$, $\nu_\star$ an eigenvaluation for $f$, and $(\pi,p,\hat{f})$ a rigidification obtained from $\nu_\star$ as in Theorem \ref{thm:rigidification}. Let $\lambda \in \nD^*$ be the non-zero eigenvalue of $df_0$.
Then $\nu_\star$ can be only a (formal) curve valuation, and:
\begin{enumerate}[(i)]
\item if $\nu_\star$ is a (non-contracted) analytic curve valuation, then $\hat{f} \cong (\lambda z, z^c w^d)$,
with $c \geq 1$ and $d \geq 1$;
\item if $\nu_\star$ is a non-analytic curve valuation, then $\hat{f} \cong (\lambda z, z^q w + P(z))$, with $q \geq 1$, and $P \in z\nC[z]$ with $\deg P \leq q$, $P \not \equiv 0$.
\end{enumerate}
\end{myprop}

\begin{mydim}
The first assertion follows from Theorem \ref{thm:fixedvaluniqueness}.

\begin{enumerate}[(i)]
\item If $\nu_\star=\nu_C$ is a (non-contracted) analytic curve valuation, then directly from \cite[Theorem 5.1]{favre-jonsson:eigenval} we have that $\mc{C}^\infty(\hat{f})=E \cup \tilde{C}$ or $\mc{C}^\infty(\hat{f})=E$, and in both cases $E$ is contracted and $\tilde{C}$ is fixed by $\hat{f}$.
We also know from Proposition \ref{prop:liftDlambda} that $\on{tr}d\hat{f}_p=\lambda \neq 0$.

In the first case $\mc{C}^\infty(\hat{f})=E \cup \tilde{C}$ is reducible: so $\hat{f}$ is of class $5$.
Hence we can choose local coordinates $(z,w)$ in $p$ such that $E=\{z=0\}$, $\tilde{C}=\{w=0\}$, and
$$
\hat{f}(z,w) = (\lambda z, z^c w^d)\mbox{,}
$$
with $c \geq 1$ and $d \geq 2$.

In the second case $\mc{C}^\infty(f)=E$ is irreducible: so $\hat{f}$ is of class $2$ or $3$, but since $E$ is contracted to $0$ by $\hat{f}$, then $\hat{f}$ is of class $2$.
Hence we can choose local coordinates $(z,w)$ such that $E=\{z=0\}$, $\tilde{C}=\{w=0\}$, and
$$
\hat{f}(z,w)=\big(\lambda z, z^q w + P(z)\big)\mbox{,}
$$
with $q \geq 1$.
Since $\tilde{C}$ is fixed, then $P \equiv 0$.

\item Let us suppose $\nu_\star = \nu_C$ a non-analytic curve valuation. We have seen in the proof of Theorem \ref{thm:rigidification} that $\mc{C}^\infty(\hat{f})=E$, and $E$ is contracted to $0$ by $\hat{f}$.
We also know from Proposition \ref{prop:liftDlambda} that $\on{tr}d\hat{f}_p=\lambda \neq 0$, so $\hat{f}$ is of class $2$ or $3$.
But only for maps in class $2$ $\hat{f}$ contracts the component $E$ in $\mc{C}^\infty(\hat{f})$.
So we are in class $2$ and we can choose local coordinates $(z,w)$ at $p$ such that $E=\{z=0\}$, and
$$
\hat{f}(z,w)=\big(\lambda z, z^q w + P(z)\big)\mbox{,}
$$
with $q \geq 1$, and $P \in z\nC[z]$ with $\deg P \leq q$. 
Since $f_\bullet^n \to \nu_\star$ in $U(p)$, no analytic curve valuation (besides $\nu_z$) is fixed by $\hat{f}$, and $P \not \equiv 0$.
\end{enumerate}
\end{mydim}

\mysubsubsect{Germs of type $(0, \nC \setminus \overline{\nD})$}

\begin{myprop}\label{prop:rigidrepel}
Let $f$ be a (dominant) holomorphic germ of type $(0, \nC \setminus \overline{\nD})$, $\nu_\star$ an eigenvaluation for $f$, and $(\pi,p,\hat{f})$ a rigidification obtained from $\nu_\star$ as in Theorem \ref{thm:rigidification}.
Let $\lambda \in \nC \setminus \overline{\nD}$ be the non-zero eigenvalue of $df_0$.
Then $\nu_\star$ can be only an analytic curve valuation, and $\hat{f} \stackrel{\on{for}}{\cong} (\lambda z, z^c w^d (1+\varepsilon z^l))$, with $c \geq 1$, $d \geq 1$, $l \geq 1$ and $\varepsilon = 0$ if $\lambda^l\neq d$, or $\varepsilon \in \{0,1\}$ if $\lambda^l=d$.
\end{myprop}

\begin{mydim}
Thanks to Theorem \ref{thm:fixedvaluniqueness}, we know that $\nu_\star$ has to be a (formal) curve valuation.

Let us suppose $\nu_\star=\nu_C$ a non-analytic curve valuation. From the proof of Theorem \ref{thm:rigidification} and Proposition \ref{prop:basinofattraction}.(i) we know that $f_\bullet^n \to \nu_C$ on a suitable open set $U=U(p)$, and hence $\hat{f}_\bullet^n \to \nu_{\tilde{C}}$ on $\mc{V}\setminus \nu_E$, where $\tilde{C}$ is the strict transform of $C$ (and it is non-analytic as well). Notice that $\nu_E$ is an analytic curve valuation if considered on the valuative tree where $\hat{f}_\bullet$ acts.
In particular, $E$ is the only analytic curve fixed by $\hat{f}$, that is in contradiction with the Stable/Unstable Manifold Theorem (see \cite[Theorem 3.1.2 and Theorem 3.1.3]{abate:hypdynsys}), since we know from Proposition \ref{prop:liftDlambda} that $\on{Spec}(d\hat{f}_p)=\{0,\lambda\}$ and $\abs{\lambda}>1$.
So $\nu_\star=\nu_C$ is a (non-contracted) analytic curve valuation.

Then the assertion on normal forms follows from Theorem \ref{thm:rigidsemiattr}.

\end{mydim}

\mysubsubsect{Germs of type $(0, \partial \nD)$}

\begin{myprop}\label{prop:rigidindif}
Let $f$ be a (dominant) holomorphic germ of type $(0, \partial \nD)$, $\nu_\star$ an eigenvaluation for $f$, and $(\pi,p,\hat{f})$ a rigidification obtained from $\nu_\star$ as in Theorem \ref{thm:rigidification}.
Let $\lambda \in \partial \nD$ be the non-zero eigenvalue of $df_0$.
Then $\nu_\star$ can be only a (formal) curve valuation, and:
\begin{enumerate}[(i)]
\item if $\lambda$ is not a root of unity, then $\hat{f} \stackrel{\on{for}}{\cong} (\lambda z, z^c w^d)$, with $c,d \geq 1$;
\item if $\lambda^r=1$ is a root of unity, then $\hat{f} \stackrel{\on{for}}{\cong} \big(\lambda z (1 + z^s + \beta z^{2s}), z^c w^d(1+\varepsilon(z^r))\big)$, where $c,d \geq 1$, $r|s$ and $\beta \in \nC$, while $\varepsilon$ is a formal power series in $z^r$, and $\varepsilon \equiv 0$ if $d \geq 2$
\end{enumerate}
\end{myprop}

\begin{mydim}
The first assertion follows from Theorem \ref{thm:fixedvaluniqueness}, while the normal forms are given by Theorem \ref{thm:rigidsemiattr}.
\end{mydim}

\mysubsect{Some remarks and examples}

\begin{myoss}\label{oss:simplerigid}
The proof of Theorem \ref{thm:rigidification} gives a general procedure to obtain a rigid germ. But in specific instances we can choose an infinitely near point lower that the one indicated.
In particular, if $\nu_\star$ is divisorial, it can happen that $U=U(p)$ can be associated to a free point $p$, and not to a satellite one.
If this is the case, we obtain a irreducible rigid germ, of class $2$ or $3$; and it has to be of class $3$, since the generalized critical set $E$ is fixed by $\hat{f}$.
So, for example, if $\hat{f}$ is still attracting, then $\hat{f}\cong (z^p, \alpha w)$, with $p \geq 2$ and $0 < \abs{\alpha} < 1$ (with $\alpha=\lambda$ if $df_0=D_\lambda$).
\end{myoss}

\begin{myes}
We present an example of the phenomenon we described in Remark \ref{oss:simplerigid}.
Set
$$
f(z,w)=(z^n+w^n, w^n)\mbox{,}
$$
with $n \geq 2$ an integer.
We easily see that $\nu_\mg{m}$ is an eigenvaluation for $f$.
We want to study the action of $\hat{f}$ on the exceptional component $E=E_0$ that arises from the single blow-up of the origin.
We can study it by checking the action of $f_\bullet$ on $\mc{E}:=\{\nu_{y-\theta x}\ |\ \theta \in \nC\} \cup \{\nu_x\}$, where we fix the corrispondence $\theta \mapsto \nu_{y-\theta x}$ between $E \cong \nP^1(\nC)$ and $\mc{E}$ (setting $\infty \mapsto \nu_x$).
Direct computations show that
$$
\hat{f}|_E : \theta \mapsto \frac{\theta^n}{1+\theta^n} \mbox{.}
$$
Now set $p=\theta \in E$ such that $\theta$ is a non-critical fixed point for $\hat{f}|_E$, i.e., such that $\theta^n + 1 = \theta^{n-1}$, and lift $f$ to a holomorphic germ $\hat{f}$ on the infinitely near point $p$.
Using the same arguments as in the proof of Theorem \ref{thm:rigidification}, we can tell that $\hat{f}$ is a rigid germ.

We show this claim by direct computations. Let us make a blow-up in $0 \in \nC^2$:
$$
\left\{
\begin{array}{l}
z=u \mbox{,} \\
w=ut \mbox{;}
\end{array}
\right.
\qquad
\left\{
\begin{array}{l}
u=z \mbox{,} \\
t=w/z \mbox{;}
\end{array}
\right.
$$
we obtain
$$
\hat{f}(u,t)=\left(u^n(1+t^n), \frac{t^n}{1+t^n}\right)\mbox{.}
$$
Choosing the local coordinates $(u, v:=t-\theta)$, we obtain
$$
\hat{f}(u,v)=	\Big(u^n\big(1+(v+\theta)^n\big), v \xi(v)\Big)\mbox{,}
$$
for a suitable invertible germ $\xi$.
In particular, $\hat{f}$ is a rigid germ, it belongs to class $3$, and (by direct computation) it is locally holomorphically conjugated to $(u,v) \mapsto (u^n, \alpha v)$, for a suitable $\alpha \neq 0$, whereas \cite[Theorem 5.1]{favre-jonsson:eigenval} would give us a germ that belongs to class $5$.
In this case, we recover the result of \cite[Theorem 5.1]{favre-jonsson:eigenval} simply by taking the lift of $g=\hat{f}$ when we blow-up the point $[0:1]\in E$, and obtaining
$$
\hat{g}(x,y)=\big(x^n y^{n-1} \chi(y),v \xi(v)\big)\mbox{,}
$$
for a suitable invertible germ $\chi$; this germ is locally holomorphically conjugated to $(x^n y^{n-1}, \alpha y)$.
\end{myes}

\begin{myoss}\label{oss:avoidedclasses}
We can apply \cite[Theorem 5.1]{favre-jonsson:eigenval}, Propositions \ref{prop:rigidattr}, \ref{prop:rigidrepel} and \ref{prop:rigidindif} even when $f$ is rigid itself: the result is that we can avoid some kind of rigid germs.
First of all, from the proof of \cite[Theorem 5.1]{favre-jonsson:eigenval} (and recalling Proposition \ref{prop:liftDlambda}), one can see that Class $7$ can be always avoided (hence Class $7$ is not ``stable under blow-ups'').
Moreover, from the proof of Theorem \ref{thm:rigidification}, we see that the germs we obtain after lifting are such that $\hat{f}_\bullet$ has always only one fixed point $\mu_\star=\pi_\bullet^{-1}\nu_\star$ of the same type of $\nu_\star$, with two exceptions: either $\nu_\star$ is divisorial, and $\mu_\star$ turns out to be an analytic curve valuation (contracted by $\pi$), or $\nu_\star$ is an irrational eigenvaluation, and in this case it can happen that $\hat{f}_\bullet=\on{id}$ on $[\nu_z, \nu_w]$.

In the first case reapplying Propositions \ref{prop:rigidattr}, \ref{prop:rigidrepel} and \ref{prop:rigidindif} we see that we obtain the same type of germ.
In the second case, we have, up to local holomorphic conjugacy, that $\hat{f}(z,w)=(z^n, w^n)$, with a suitable $n \geq 2$. Then all valuations on $[\nu_z, \nu_w]$ are eigenvaluations, and reapplying \cite[Theorem 5.1]{favre-jonsson:eigenval}, we obtain a rigid germ that belongs to a different class. In particular, making a single blow-up on the origin, and considering the germ at $[1:1]$, we obtain a germ of the form $(nz(1+h(z)), w^n)$ for a suitable holomorphic map $h$ such that $h(0)=0$, that is (by direct computation) holomorphically conjugated to $(nz, w^n)$.
\end{myoss}

\begin{myes}\label{es:normalformchanges}
Reapplying \cite[Theorem 5.1]{favre-jonsson:eigenval}, as just seen in the last remark, we usually obtain the same normal form type. But there are cases where the normal form can change (staying rigid).

Consider for example the rigid germ $f(z,w)=(w^2, z^3)$. Then the only eigenvaluation $\nu_\star$ is the monomial valuation on $(z,w)$, such that $\nu_\star(x)=1$ and $\nu_\star(w)=\sqrt{3/2}$.
Then an infinitely near point $p$ that works in \cite[Theorem 5.1]{favre-jonsson:eigenval} can be obtained after three blow-ups:
the first at $0$ (and we obtain $E_0$), the second at $[1:0] \in E_0$ (and we obtain $E_1$), the third at $[0:1]$ (and we obtain $E_2$). We can choose $p=[0,1] \in E_2$, and the lift we obtain is $\hat{f}(z,w)=(w^6, z)$.
\end{myes}

\nocite{zariski-samuel:commalg1}
\nocite{hubbard-papadopol:supfixpnt}

\bibliographystyle{alpha}
\bibliography{biblio}

\end{document}